\magnification=1200     \hfuzz=1pt
\hsize=125mm \vsize=187mm \hoffset=4mm \voffset=10mm

\pretolerance=500 \tolerance=1000 \brokenpenalty=5000
\magnification=1200
\hsize=125mm \vsize=187mm \hoffset=4mm \voffset=10mm
\parindent=0mm

\catcode`\$=3
\catcode`\@=11
\catcode`\;=\active
\def;{\relax\ifhmode\ifdim\lastskip>\z@\unskip\fi\kern.2em\fi\string;}
\catcode`\:=\active
\def:{\relax\ifhmode\ifdim\lastskip>\z@\unskip\fi\penalty\@M\ \fi\string:}
\catcode`\!=\active
\def!{\relax\ifhmode\ifdim\lastskip>\z@\unskip\fi\kern.2em\fi\string!}
\catcode`\?=\active
\def?{\relax\ifhmode\ifdim\lastskip>\z@\unskip\fi\kern.2em\fi\string?}
\catcode`\@=12

\font\twbf=cmbx12
\def\vk#1{\vskip #1mm}
\parskip=0pt plus 1pt minus 1pt
\def\parag{\vk{2}}
\def\qed{\hfill\kern 6pt
\lower 2pt\hbox{\vrule\vbox to10pt{\hrule width 4pt \vfil\hrule}\vrule}\par}

\def\cf{{\sl cf.\/}\ }
\def\op#1{\mathop{\rm #1}}
\def\e#1{\op{e}\nolimits^{#1}}
\def\Lim{\op{Lim}}
\def\Sup{\op{Sup}}
\def\R{{\rm I\! R}}
\def\Frac#1#2{\displaystyle{{#1}\over {#2}}}
\let\eps=\varepsilon
\let\ph=\varphi
\let\cv=\rightarrow
\let\imp=\Rightarrow
\def\parno{\count101}
\def\bibno{\count102}
\long\def\bib#1#2#3{\item{[\the\bibno]}{\underbar{#1}}~{\sl #2}~{#3\par}
\advance\bibno by1}
\def\sect#1{{\bf\the\parno\ #1:~}\global\advance\parno by1}

\def\thm{\sect{Theorem}\sl}
\def\cor{\sect{Corollary}\sl}
\def\defi{\sect{Definition}\sl}
\def\rem{\sect{Remark}}
\def\rems{\sect{Remarks}}
\def\prop{\sect{Proposition}\sl}
\def\exple{\sect{Example}\sl}
\def\dm{\rm Proof~:~}
\def\cC{{\cal C}}
\let\adh=\overline
\let\lt=\left
\let\rt=\right
\let\otim=\otimes
\def\cA{{\cal A}}
\def\cM{{\cal M}}
\def\cU{{\cal U}}\hfill
\def\C{\rm\;{}^{{}_|}\!\!\! C\,}

\bibno=1
\parno=1

\def\,{\mkern2mu}
\def\;{\mkern4mu}
\hfill {4 juin 2007}\parag

\hfil {\twbf A NON-COMMUTATIVE SEWING LEMMA}     \parag  
\hfil { Denis Feyel\par 
\hfil D\'epartement de Math\'ematiques, Universit\'e d'Evry-Val d'Essonne\par 
\hfil Boulevard Francois Mitterrand, 91025 Evry cedex, France \par 
\hfil Denis.Feyel@univ-evry.fr}   \parag  
\hfil {Arnaud de La Pradelle\par 
\hfil Laboratoire d'Analyse Fonctionnelle, Universit\'e Paris VI\par 
\hfil Tour 46-0, 4 place Jussieu, 75052 Paris, France\par 
\hfil Universit\'e Pierre-et-Marie Curie\par 
\hfil adelapradelle@free.fr}   \parag  
\hfil {Gabriel Mokobodzki\par 
\hfil Laboratoire d'Analyse Fonctionnelle, Universit\'e Paris VI\par 
\hfil Tour 46-0, 4 place Jussieu, 75052 Paris, France\par 
\hfil Universit\'e Pierre-et-Marie Curie}   \parag     \parag


{\bf Abstract} A non-commutative version of the sewing lemma [1] is proved,
with some applications.   \parag     \parag

{\bf Introduction}   \parag  
In a preceding paper [1] we proved a sewing lemma which was a key result for the
study of H\"older continuous functions. In this paper we give a
non-commutative version of this lemma.  \vk{1}
In the first section we recall the commutative version, and give some
applications (Young integral and stochastic integral).  \vk{1}
In the second section we prove the non-commutative version. This last result
has interesting applications: an extension of the so-called integral
product, a simple case of the semigroup Trotter type formula, and a sharpening
of the Lyons theorem about multiplicative functionals [3,4,5].  \vk{1}
Note that we replaced the H\"older modulus of continuity $t^\alpha $ by a more
general modulus $V(t)$.   \parag  
This paper was elaborated with the regretted G. Mokobodzki. The writing has
only been done after his death.   \parag     \parag     \parag

{\bf I. The additive sewing lemma}   \parag

\defi{} We say that a function $V(t)$ defined on $[0,T[$ is a control
function if it is non decreasing, $V(0)=0$ and $\sum  _{n\ge 1}V(1/n)<\infty $.   \parag  \rm

As easily seen, this is equivalent to the property
$$\adh V(t)=\sum  _{n\ge 0}2^nV(t.2^{-n})<\infty $$
for every $t\ge 0$. For example, $t^\alpha $ and $t/(\log{t^{-1}})^\alpha $ with $\alpha >1$ are
control functions.   \parag
Observe that we have
$$\adh V(t)=V(t)+\cdots+2^nV(t.2^{-n})+2^{n+1}\adh V(t.2^{-n-1})$$
from which follows that $\displaystyle\Lim_{t\cv 0}\adh V(t)/t=0$.   \parag     \parag

\thm{} Consider a continuous function $\mu (a,b)$ defined for $0\le a\le b<T$ satisfying the relation
$$|\mu (a,b)-\mu (a,c)-\mu (c,b)|\le V(b-a)$$
for every $c\in [a,b]$,
where $V$ is a control function. Then there exists a unique function $\ph (t)$
on $[0,T[$, up to an additive constant, such that
$$|\ph (b)-\ph (a)-\mu (a,b)|\le \adh V(b-a)$$
\dm Put $\mu '(a,b)=\mu (a,c)+\mu (c,b)$ for $c=(a+b)/2$, and $\mu ^{(n+1)}=\mu ^{(n)}{}'$.
We easily get for $n\ge 0$
$$|\mu ^{(n)}(a,b)-\mu ^{(n+1)}(a,b)|\le 2^nV(2^{-n}|b-a|)$$
so that the series $\sum  _{n\ge 0}|\mu ^{(n)}(a,b)-\mu ^{(n+1)}(a,b)|\le \adh V(b-a)$ converges, and the
sequence $\mu ^{(n)}(a,b)$ converges to a limit $u(a,b)$. For $c=(a+b)/2$ we
have $\mu ^{(n+1)}(a,b)=\mu ^{(n)}(a,c)+\mu ^{(n)}(c,b)$ which implies
$$u(a,b)=u(a,c)+u(c,b)$$
We say that $u$ is midpoint-additive.  \vk{1}
Now, we prove that $u$ is the unique midpoint-additive function with the inequality
$|u(a,b)-\mu (a,b)|\le \op{Cst}\adh V(b-a)$. Indeed if we have another one $v$, we get
$$|v(a,b)-u(a,b)|\le K.\adh V(b-a)$$
and by induction $|v(a,b)-u(a,b)|\le 2^nK.\adh V[2^{-n}(b-a)]$ which vanishes as $n\cv \infty $
as mentioned above. Let $k$ be an integer $k\ge 3$, and take the function
$$w(a,b)=\sum  _{i=0}^{k-1}u(t_i,t_{i+1})$$
with $t_i=a+i.(b-a)/k$. It follows that $w$ also is midpoint-additive, and
satisfies
$$|w(a,b)-\mu (a,b)|\le \hbox{Cst}_k\,\adh V(b-a)$$
hence we have $w=u$, that is $u$ is in fact rationally-additive. As $\mu $ is
continuous, then so also is $u$, as the defining series converges uniformly
for $0\le a\le b<T$. Then $u$ is additive, and it suffices to put
$\ph (t)=u(0,t)$.\qed    \parag     \parag

\prop{(Riemann sums)} Let $\sigma =\{t_i\}$ some finite subdivision of $[a,b]$.
Put $\delta =\Sup_i|t_{i+1}-t_i|$. Then
$$\Lim_{\delta \cv 0}\;\sum  _i\mu (t_i,t_{i+1})=\ph (b)-\ph (a)$$
\dm We have
$$\ph (b)-\ph (a)-\sum  _i\mu (t_i,t_{i+1})=\sum  _i[\ph (t_{i+1})-\ph (t_i)-\mu (t_i,t_{i+1})]$$
$$\lt|\ph (b)-\ph (a)-\sum  _i\mu (t_i,t_{i+1})\rt|\le \sum  _i\adh V(t_{i+1}-t_i))\le \eps \sum  _i(t_{i+1}-t_i)=(b-a)\eps $$
since $\adh V(\delta )/\delta \le \eps $ as $\delta \cv 0$.   \parag     \parag

\rems{} a) In fact the result holds even in the case of discontinuous $\mu $, as
proved in the appendix.   \vk{1}
b) The result obviously extends to Banach spaces valued functions $\mu $.   \parag
In the case $V(t)=t^\alpha $ with $\alpha >1$, we get $\adh V(t)=\Frac{2^\alpha \,t^\alpha }{2^\alpha -2}$.   \parag     \parag  

\exple The Young integral   \parag
Take $V(t)=t^{2\alpha }$ with $\alpha >1/2$. If $x$ and $y$ are two $\alpha $-H\"older
continuous functions on $[0,1]$, put
$$\mu (a,b)=x_a(y_b-y_a)$$
We get
$$\mu (a,b)-\mu (a,c)-\mu (c,b)=-(x_c-x_a)(y_b-y_c)$$
so that
$$|\mu (a,b)-\mu (a,c)-\mu (c,b)|\le \|x\|_\alpha \|y\|_\alpha |b-a|^{2\alpha }$$
where $\|x\|_\alpha $ is the norm in the space $\cC^\alpha $. Let $\ph $ be the function of
theorem 2, put
$$\int _a^bx_t\,dy_t=\ph (b)-\ph (a)$$
This is a Young integral (\cf also [7]).   \parag  
\rem  We could take $x\in \cC^\alpha $, $y\in \cC^\beta $ with $\alpha +\beta >1$.   \parag     \parag

\exple The stochastic integral   \parag
Let $X_t$ be the standard $\R^m$-valued Brownian motion. As is well known,
$t\cv X_t$ is $\cC^{1/2}$ with values in $L^2$. Let $f$ be a tensor-valued
$\cC^2$-function with bounded derivatives on $\R^m$. Put
$$\mu (a,b)=f(X_a)\otim (X_b-X_a)+\nabla{f}(X_a)\otim \int _a^b(X_t-X_a)\otim dX_t$$
where the last integral is taken in the Ito or in the Stratonovitch sense.
By straightforward computations, we get
$$N_2[\mu (a,b)-\mu (a,c)-\mu (c,b)]\le K.\|\nabla{f}\|_{\cC^1}|b-a|^{3/2}$$
By the additive sewing lemma, there exists a unique $L^2$-valued function
$\ph (t)$ such that $N_2[\ph (b)-\ph (a)-\mu (a,b)]\le \op{Cst}.|b-a|^{3/2}$ (the control
function is $V(t)=t^{3/2}$). It is easily seen that
$$\ph (b)-\ph (a)=\int _a^bf(X_t)\otim dX_t$$
in the Ito or in the Stratonovitch sense.   \parag  
Observe that as the stochastic integral $\displaystyle\int _a^bX_t\otim dX_t$ has
$\cC^\alpha $-trajectories almost surely for $1/3<\alpha <1/2$, analoguous computations
as above yield $\cC^\alpha $-trajectories for $\displaystyle\int _a^bf(X_t)\otim dX_t$ on the same
set of paths as $\displaystyle\int _a^bX_t\otim dX_t$.   \parag     \parag

\rem{} For the FBM with $\alpha >1/4$, the reader is referred to our previous
paper [1].   \parag     \parag  

{\bf II. The multiplicative sewing lemma}   \parag  
Here we need a strong notion of control function   \parag  

\defi{} We say that a function $V(t)$ defined on $[0,T[$ is a strong control
function if it is a control function and there exists a $\theta >2$ such that
for every $t$
$$\adh V(t)=\sum  _{n\ge 0}\theta ^nV(t.2^{-n})<\infty $$\rm
We consider an associative monoide $\cM$ with a unit element $I$, and we
assume that $\cM$ is complete under a distance $d$ satisfying
$$d(xz,yz)\le |z|\,d(x,y),~~~~~~d(zx,zy)\le |z|\,d(x,y)$$
for every $x,y,z\in \cM$, where $z\cv |z|$ is a Lipschitz function on $\cM$ with
$|I|=1$.   \parag  
Let $\mu (a,b)$ be an $\cM$-valued function defined for $0\le a\le b<T$. We assume
that $\mu $ is continuous, that $\mu (a,a)=I$ for every $a$, and that for every
$a\le c\le b$ we have
$$d(\mu (a,b),\mu (a,c)\mu (c,b))\le V(b-a)\leqno  (1)$$
We say that an $\cM$-valued $u(a,b)$ is multiplicative if
$u(a,b)=u(a,c)u(c,b)$ for every $a\le c\le b$.   \parag  
\thm{} There exists a unique multiplicative function $u$ such that
$d(\mu (a,b),u(a,b))\le \op{Cst}\adh V(b-a)$ for every $a\le b$.   \parag  
\dm Put $\mu _0=\mu $ and by induction
$$\mu _{n+1}(a,b)=\mu _n(a,c)\mu _n(c,b)~~~~\hbox{where}~~~~c=(a+b)/2$$
$$h_n(t)=\Sup_{b-a\le t}|\mu _n(a,b)|,~~~~~~~~
U_n(t)=\Sup_{b-a\le t}d(\mu _{n+1}(a,b),\mu _n(a,b))$$
The functions $h_n$ and $U_n$ continuous non decreasing with
$h_n(0)=1$ and $U_n(0)=0$. Let $\kappa $ be the Lipschitz constant of $z\cv |z|$.
One has
$$h_{n+1}(t)\le h_n(t)+\kappa U_n(t)\le h_0(t)+\kappa U_0(t)+\cdots+\kappa U_n(t)$$
$$U_{n+1}(t)\le 2h_{n+1}(t/2)U_n(t/2)\leqno  (2)$$

Let $\tau >0$ be such that $h_0(\tau )+\kappa \adh V(\tau )\le \theta /2$. Assume that
$U_i(t)\le \theta ^iV(t/2^i)$ for $t\le \tau $ and $i\le n$. One has $h_{n+1}(t)\le \theta /2$, then
$$U_{n+1}(t)\le \theta U_n(t/2)\le \theta ^{n+1}V(t/2^{n+1})$$
for $t\le \tau $ and every $n$ by induction.  \vk{1}
Hence for $t\le \tau $ the series $U_n(t)$ converges, so that the sequence
$h_n(\tau )$ is bounded. By inequality (2) the series $U_n(2\tau )$ converges, and
the sequence $h_n(2\tau )$ is bounded. From one step to the other we see that
the sequence $h_n$ is locally bounded, and that the series $U_n$ converges
locally uniformly on $[0,T[$. It follows that the sequence $\mu _n(a,b)$
converges locally uniformly to a continuous function $u(a,b)$ which is
midpoint-multiplicative, that is $u(a,b)=u(a,c)u(c,b)$ for $c=(a+b)/2$.
One has $d(u,\mu )\le \op{Cst}\adh V$.  \vk{1}
Next we prove the unicity of $u$. Let $v$ be a function with the same
properties as $u$. Put $K(t)=\displaystyle\Sup_{b-a\le t}[u(a,b),v(a,b)]$. Let $\tau _1>0$ be
such that $K(\tau _1)\le \theta /2$. One has
$d(u(a,b),v(a,b))\le k\,\adh V(b-a)$ with some constant $k$, then
$d(u(a,b),v(a,b))\le 2K(t/2)k\adh V(t/2)\le k\,\theta \adh V(t/2)$ for $b-a\le t\le \tau _1$, and by induction
$d(u(a,b),v(a,b))\le k\,\theta ^n\adh V(t/2^n)$. It follows that $u(a,b)=v(a,b)$
for $b-a\le \tau _1$. This equality extends to every $b-a$ by midpoint-multiplicativity.  \vk{1}
Finally we prove that $u$ is multiplicative. We argue as in the additive
case, and we put for an integer $m$
$$w(a,b)=\prod _{i=0}^{m-1}u(t_i,t_{i+1})$$
where $t_i=a+i.(b-a)/m$. For simplicity we limit ourselves to the case $m=3$,
that is
$$w(a,b)=u(a,c')u(c',c'')u(c'',b)$$
with $c'=a+(b-a)/3$, $c''=a+2(b-a)/3$. Observe that $w$ is obviously
midpoint-multiplicative. Take $a\le b\le T_0<T$, we get successively
with a constant $k$ which can be changed from line to line
$$\eqalign{d(w(a,b),\mu (a,b))\le &k\adh V(b-a)+d(w(a,b),\mu (a,c')\mu (c',b))\cr
\le k\adh V(b-a)&+d(u(a,c')u(c',c'')u(c'',b),u(a,c')\mu (c',b))\cr
&+d(\mu (u(a,c')\mu (c',b),\mu (a,c')\mu (c',b))\cr
\le k\adh V(b-a)&+kd(u(c',c'')u(c'',b),\mu (c',b))+kd(\mu (a,c'),\mu (a,c'))\cr
\le k\adh V(b-a)&+kd(u(c',c'')u(c'',b),\mu (c',b))\cr
\le k\adh V(b-a)&+kd(u(c',c'')u(c'',b),u(c',c'')\mu (c'',b))\cr
&+kd(u(c',c'')\mu (c'',b),\mu (c',c'')\mu (c'',b))\cr
\le k\adh V(b-a)&+kd(u(c'',b),\mu (c'',b))+kd(u(c',c''),\mu (c',c''))\cr
\le k\adh V(b-a)&}$$
By the second step of the proof, we get $w=u$. The same proof extends to
every $m$, so that $u$ is in fact rationally multiplicative. As $u$ is
continuous, it is multiplicative.   \qed   \parag     \parag  

\exple The integral product   \parag
Let $t\cv A_t$ a $\cC^\alpha $ function with values in a Banach algebra $\cA$ with a
unit $I$. Put $A_{ab}=A_b-A_a$ and
$$\mu (a,b)=I+A_{ab}$$
We get
$$\mu (a,b)-\mu (a,c)\mu (c,b)=-A_{ac}A_{cb}$$
Suppose that $\alpha >1/2$, then the multiplicative sewing lemma applies with the
obvious distance, and there exists a unique multiplicative function $u(a,b)$
with values in $\cA$ such that
$$|u(a,b)-\mu (a,b)|\le \op{Cst}|b-a|^{2\alpha }$$
We get the same $u(a,b)$ by taking $\mu (a,b)=\e{A_{ab}}$.
A good notation for $u(a,b)$ is
$$u(a,b)=\prod _a^b(I+dA_t)=\prod _a^b\e{dA_t}$$
\thm{} Put $H_t=u(0,t)$. Then this is the solution of the EDO
$$H_t=I+\int _0^tH_s\,dA_s$$
\dm We have only to verify that
$|u(0,b)-u(0,a)-u(0,a)A_{ab}|\le \op{Cst}|b-a|^{2\alpha }$. The first member is worth
$$u(0,a)[u(a,b)-I-A_{ab}]=u(0,a)[u(a,b)-\mu (a,b)]$$
so that we are done.   \parag     \parag

\exple A Trotter type formula   \parag
Let $t\cv A_t$ and $t\cv B_t$ as in the previous paragraph, and put
$$\mu (a,b)=[I+A_{ab}][I+B_{ab}]$$
It is straightforward to verify the good inequality
$$|\mu (a,b)-\mu (a,c)\mu (c,b)|\le \op{Cst}|b-a|^{2\alpha }$$
so that we get a multiplicative $u(a,b)$ such that
$|u(a,b)-\mu (a,b)|\le \op{Cst}|b-a|^{2\alpha }$ or equivalently
$$|u(a,b)-I-A_{ab}-B_{ab}|\le \op{Cst}|b-a|^{2\alpha }$$
$$u(a,b)=\prod _a^b(I+dA_t+dB_t)=\prod _a^b\e{dA_t}\e{dB_t}$$
We then have
$$\e{A_{ab}+B_{ab}}=\Lim_{n\cv \infty }\prod _{i=1}^n\e{A_{t_it_{i+1}}}\e{B_{t_it_{i+1}}}$$
for $t_{i+1}-t_i=(b-a)/2^n$.   \parag  
Particularly we can take $A_t=tA$ and $B_t=tB$ with $\alpha =1$, this yields
$$\e{A+B}=\Lim_{n\cv \infty }\prod _{i=1}^n\e{A/2^n}\e{B/2^n}$$
   \parag     \parag

\exple Extending the Lyons theorem   \parag
Let $\cA$ be a Banach algebra with a unit $I$.
Take $\mu (a,b)$ of the form
$$\mu (a,b)=\sum  _{k=0}^n\lambda ^kA_{ab}^{(k)}$$
where $A_{ab}^{(k)}\in \cA$, $\lambda $ is a real parameter. We have
$$\mu (a,c)\mu (c,b)=\sum  _{k=0}^n\lambda ^kB_{acb}^{(k)}+\sum  _{k=n+1}^{2n}\lambda ^kC_{acb}^{(k)}$$
Following [5], we suppose the algebraic hypothesis for $k\le n$
$$A^{(k)}_{ab}=\sum  _{i=0}^kA^{(i)}_{ac}A^{(k-i)}_{cb}\leqno  (3)$$
that is
$$\mu (a,c)\mu (c,b)=\mu (a,b)+\sum  _{k=n+1}^{2n}\lambda ^kC_{acb}^{(k)}$$

\thm{} Under the condition (3) and the inequality
$$|A^{(k)}_{ab}|\le M|b-a|^{k\alpha }$$
for every $k\le n$, where $\alpha >1/(n+1)$, there exists a unique multiplicative
function $u(a,b)$ such that
$$|u(a,b)-\mu (a,b)|\le \op{Cst}|b-a|^{(n+1)\alpha }$$
Moreover we have
$$u(a,b)=\sum  _{k=0}^n\lambda ^kA_{ab}^{(k)}+\sum  _{n+1}^\infty \lambda ^kB_{ab}^{(k)}\leqno  (4)$$
where the series is normally convergent for every $\lambda $.   \parag  
\dm The only problem is to prove formula (4), that is to prove that $u$ is
the sum of its Taylor expansion with respect to $\lambda $. In the case where $\cA$
is a complex Banach algebra, the proof of the multiplicative sewing lemma yields a
sequence of holomorphic functions which converges uniformly with respect to
$\lambda $ in every compact set of $\C$. Hence $u(a,b)$ in holomorphic in $\lambda \in \C$.
If $\cA$ is only a real Banach algebra, we get a sequence of holomorphic
functions with values in the complexified Banach space of $\cA$, and the
result follows.  It remains to observe that the $n+1$ first terms of the
Taylor expansion are the same for every function of the sequence converging to
$u$.   \parag     \parag  
{\bf Application to the Lyons theorem:} Let $E$ be a Banach space. Denote $E^n=E^{\otim n}$. For every
$k\le n$, let $(a,b)\cv X^{(k)}_{ab}$ an $E^k$-valued function such that
$$X^{(k)}_{ab}=\sum  _{i=0}^{k}X^{(i)}_{ac}\otim X^{(k-i)}_{cb}$$
for $a\le c\le b$. Suppose that every $E^n$ has a cross-norm such that
$$\|u\otim v\|_{n+m}\le \|u\|_n\|v\|_m$$
for every $u\in E^n$, $v\in E^m$. Suppose that $\alpha >1/(n+1)$, and that we have for $k\le n$
$$\|X^{(k)}_{ab}\|_k\le M.|b-a|^{k\alpha }$$
Let $\cA$ be the completed tensor algebra under the norm
$$\|t\|=\sum  _{n\ge 0}\|t_n\|_n$$
This is a Banach algebra. The previous theorem applies, so that
there exists a unique $(a,b)\cv Y^{(k)}_{ab}$ for every $k$ such that
$Y^{(k)}=X^{(k)}$ for $k\le n$,
$$Y^{(k)}_{ab}=\sum  _{i=0}^{k}Y^{(i)}_{ac}\otim Y^{(k-i)}_{cb}$$
for every $k$ and every $a\le c\le b$, and
$$\sum  _{k\ge n+1}\|Y^{(k)}_{ab}\|_k\le \op{Cst}|b-a|^{(n+1)\alpha }$$   \parag  
\rem This theorem sharpens the theorem 3.2.1 of [5].   \parag     \parag     \parag

{\bf Some estimations}   \parag  
We return to formula (4) of theorem 9
$$u(a,b)=\sum  _{k=0}^N\lambda ^kA_{ab}^{(k)}+\sum  _{N+1}^\infty \lambda ^kB_{ab}^{(k)}$$
for $N=\op{Ent}(1/\alpha )$,
and we put $B_{ab}^{(k)}=A_{ab}^{(k)}$ for simplification, so that we have
$$u(a,b)=\sum  _{k=0}^\infty \lambda ^kA_{ab}^{(k)}$$
There exist best constants $K_n$ such that
$|A^{(n)}_{ab}|\le K_n|b-a|^{n\alpha }$. We have
$$A^{(n+1)}_{ab}=A^{(n+1)}_{ac}+A^{(n+1)}_{cb}+\sum  _{k=1}^nA^{(k)}_{ac}A^{(n-k+1)}_{cb}$$
By taking $c=(a+b)/2$ we get
$$|A^{(n+1)}_{ab}|\le 2^{-(n+1)\alpha }\lt[2K_{n+1}+\sum  _{k=1}^nK_kK_{n-k+1}\rt]|b-a|^{(n+1)\alpha }$$
and then
$$(2^{(n+1)\alpha }-2)K_{n+1}\le \sum  _{k=1}^nK_kK_{n-k+1}$$
Let $0<\beta <\alpha $, and introduce the entire function
$$\e{}(x)=\e{}_\beta (x)=\sum  _{n\ge 0}\Frac{x^n}{n!^\beta }~~~~~~\imp ~~~~~~\e{}(x)^2=\sum  _{n\ge 0}E_{n,\beta }\Frac{x^n}{n!^\beta }$$
where
$$E_{n,\beta }=\sum  _{k=0}^n\lt[C_n^k\rt]^\beta \le 2^{n\beta }(n+1)$$
There exist $c\ge 0$ and $x>0$ such that for $1\le m\le N$
$$K_m\le c.x^m/m!^\beta \leqno  (5)$$
Hence we have for $n\ge N$
$$(2^{(n+1)\alpha }-2)K_{n+1}\le c^2x^{n+1}\sum  _{k=1}^n(k!)^{-\beta }(n-k+1)^{-\beta }\le c^2x^{n+1}[(n+1)!]^{-\beta }A_{n+1,\beta }$$
In order that (5) holds for every $n$, it suffices that
$${1\over{c}}\ge \Sup_{n>N}\Frac{A_{n+1,\beta }}{2^{(n+1)\alpha }-2}$$
which is possible since the fraction in the second hand member shrinks to 0
as $n\cv \infty $.   \parag  
\cor{} Put $c'=\op{Max}(c,1)$, we have
$$|u(a,b)|\le c'\e{}_\beta (|\lambda |x|b-a|^\alpha )$$   \parag     \parag  \rm

\rems{} a) Note that for $\alpha =1$ one can take $\beta =\alpha =1$ so that we recover the
classical inequality.  \vk{1}
b) For $\beta <1$, the function $\e{}_\beta (x)$ increases faster than the exponential
function (\cf Schwartz [6] for $\beta =1/2$).  \vk{1}
c) there are some analoguous computations in Gubinelli [2].   \parag  
        \parag     \parag  
{\bf Appendix: the discontinuous case}   \parag  
As announced in Remark 4a), we extend the additive sewing lemma in the case
where $\mu $ is discontinuous. We go back to the proof of the lemma: we get a
unique function $u(a,b)$ which is rationally additive and such that
$|u(a,b)-\mu (a,b)|\le \op{Cst}\adh V(b-a)$. Put
$$v_n(a,b)=u(a_n,b_n)-u(a_n,a)+u(b_n,b)$$
where $a_n\le a$ and $b_n\le b$ are the classical dyadic approximations of $a$ and
$b$. It is straightforward to verify that $v_n$ is additive for every
$a\le c\le b$. Besides, we have
$$|v_n(a,b)-u(a,b)|\le 2\adh V(b-a)+2\adh V(b_n-a_n)+\adh V(a-a_n)+\adh V(b-b_n)$$
so that the sequence $v_n(a,b)-u(a,b)$ is bounded. Let $v(a,b)$ be the limit of
$v_n(a,b)$ according to an ultrafilter $\cU\cv \infty $. We first have
$v(a,b)=v(a,c)+v(c,b)$ for every $a\le c\le b$. Then we get
$$|v(a,b)-\mu (a,b)|\le 3\adh V(b-a)+2\Lim_{\cU}\adh V(b_n-a_n)\le 5\adh V(2(b-a))$$
As $V(2t)$ is also a control function for $\mu $, $v$ is the unique additive
function such that $|v(a,b)-\mu (a,b)|\le 5\adh V(2(b-a))$, which implies that $v=u$.
Hence $u$ is completely additive. \qed
   \parag  

Here we point out the important fact that the result also holds if $\mu $
takes values in a Banach space $B$. Indeed, the proof is exactly the same,
the last limit according to $\cU$ must be taken in the bidual $B''$ with the
topology $\sigma (B'',B')$.   \parag  

   \parag     \parag

\hfil {\bf References}   \parag

\bib{D. Feyel, A. de La Pradelle.}{Curvilinear Integrals along Enriched
Paths.}{\par Electronic J. of Prob. 34, 860-892, (2006).}\vskip1mm
\bib{M. Gubinelli.}{Controlling rough paths.}{\par
J. Func. Anal., 216, pp. 86-140, (2004).}\vskip1mm
\bib{T.J. Lyons.}{Differential equations driven by rough signals.}{\par
Rev. Math. Iberoamer. {\bf14}, 215-310, (1998).}  \vskip1mm
\bib{T.J. Lyons, Z. Qian.}{Calculus for multiplicative functionals, Ito's
formula and differential equations.}{\par Ito's stochastic calculus and
Probability theory, 233-250, Springer, Tokyo, (1996).}  \vskip1mm
\bib{T.J. Lyons, Z. Qian.}{System Control and Rough Paths.}
{\par Oxford Science Publications, (2002).}  \vskip1mm
\bib{L. Schwartz.}{La convergence de la s\'erie de Picard pour les EDS.}{\par
S\'em. prob. Strasbourg, t.23, 343-354, (1989).} \vskip1mm
\bib{L.C. Young.}{An inequality of H\"older type, connected with Stieltjes
integration.}{~~Acta Math. 67, 251-282 (1936).}  \vskip1mm

\bye
\bye